# Nonlinear dynamics of mechanical systems with friction contacts: coupled static and dynamic Multi-Harmonic Balance Method and multiple solutions.

Stefano Zucca[a]

[a]Politecnico di Torino, Department of Mechanical and Aerospace Engineering,
Corso Duca degli Abruzzi, 24, 10139, Torino, Italy
tel. +39 0110906933, fax +39 0110906999
email: stefano.zucca@polito.it

Christian Maria Firrone[b]
(Corresponding author)
[b]Politecnico di Torino, Department of Mechanical and Aerospace Engineering,
Corso Duca degli Abruzzi, 24, 10139, Torino, Italy
tel. +39 0110906912, fax +39 0110906999
email: christian.firrone@polito.it

**Abstract**

Real applications in structural mechanics where the dynamic behaviour is linear are rare. Usually, structures are made of components assembled together by means of joints whose behaviour maybe highly nonlinear. Depending on the amount of excitation, joints can dramatically change the dynamic behaviour of the whole system, and the modelling of this type of constraint is therefore crucial for a correct prediction of the amount of vibration.

The solution of the nonlinear equilibrium equations by means of the Harmonic Balance Method (HBM) is widely accepted as an effective approach to calculate the steady-state forced response in the frequency domain, in spite of Direct Time Integration (DTI). The state-of-the-art contact element used to model the friction forces at the joint interfaces is a node-to-node contact element, where the local contact compliance is modelled by means of linear springs and Coulomb's law is used to govern the friction phenomena. In the literature, when the HBM is applied to vibrating systems with joint interfaces and the state-of-the-art contact model is used, an uncoupled approach is mostly employed: the static governing equations are solved in advance to compute the pre-stress effects and then the dynamic governing equations are solved to predict the vibration amplitude of the system. As a result, the HBM steady-state solution may lead to a poor correlation with the DTI solution, where static and dynamic loads are accounted for simultaneously..

In this paper, the HBM performances are investigated by comparing the uncoupled approach to a fully coupled static/dynamic approach. In order to highlight the main differences between the two approaches, a lumped parameter system, characterized by a single friction contact, is considered in order to show the different levels of accuracy that the proposed approaches can provide for different configurations.

# 1    Introduction

Joint mechanics represents one of the most interesting issues that designers have to consider during the design of complex systems. Usually machines are the result of an assembly of many parts, and when the connected elements are coupled together through removable links like bolting, shrink fits or unilateral constraints, the interactions between two or more bodies in contact may affect the vibratory response during operations. Although linear systems are currently well predicted by simulation tools based on the Finite Element Method computations, provided that suitable model updating is performed, the prediction of the nonlinear dependence of the dynamic behaviour of structures from friction contacts remains an open issue. This is sometimes converted to additional linear constraints to the model in terms of additional stiffness or additional modal damping in order to take into account a tight or loose configurations of the contact. The conversion is mainly based on experience and empirical approaches. Different efforts have been made by several authors ([1]-[5]) to take into account the actual stick/slip/lift-off states of the contact in order i) to determine the real amount of change in the resonant frequency (stiffness contribution) and in the resonant amplitude (damping contribution) due to the friction contact and ii) to optimize the contact modelling in order to reduce the calculation time of the forced response, which must be performed iteratively due to the nonlinear nature of the phenomenon. In particular in turbo machinery design, existing joints may be optimized in order to exploit the damping contribution in order to limit the structural vibrations thanks to the dissipated energy at the contact (blade root joints [6],[7], shrouds and snubber [8]-[10]) or dampers can be purposely added to the system (underplatform dampers [11]-[14] and ring dampers [15],[16]).

The state-of-the-art contact element used to model the contact forces at the joint interfaces is a node-to-node contact element ([5]), where the local contact compliance is modelled by linear springs and Coulomb's friction law is used to model the local slip conditions of the joint interfaces.

When this contact element is used to compute the forced response of highly detailed FE models, the Direct Time Integration (DTI) is unpractical due to the large calculation time needed to find the steady-state solution, and therefore the Harmonic Balance Method (HBM) is used to approximate the periodical quantities as a Fourier series ([17]) and to obtain a set of non-linear algebraic equations in the frequency domain.

A common practice is to solve in advance the static governing equations of the non-linear system by applying the static loads (assembly preloads, centrifugal force) to compute the distribution of static normal loads that act on the joint interface, with a so-called pre-stress analysis. These static loads are then used as input parameters in the dynamic analysis of the system where only the dynamic loads are included. This two-step approach can be referred to as an 'uncoupled approach', where the static equilibrium is not influenced by the dynamic equilibrium.

In the present paper, the authors demonstrate with a simple lumped spring-mass system, that, if the uncoupled approach is used, a wide range of multiple solutions ([18]) is possible, because multiple static configurations may exist. Furthermore, the uncoupled approach may sometimes give erroneous results since the static contact loads computed in advance can differ from the actual static loads computed with the coupled approach. As a result of the inaccurately computed static load, also the prediction of the vibration amplitude of the system is inaccurate. In particular, it is demonstrated that when the contact enters the slip state, the static contact loads, computed with the pre-stress analysis, are modified by the slip phenomena, and, at the same time, the modified static loads affect the amount of damping produced by the friction forces. For this reason the authors developed an improved strategy to take into account this phenomenon (the mutual dependence of the static and dynamic contact loads) and to reduce the non-uniqueness of the solution to one output only [19]. The improved strategy was presented first in [20] where a typical industrial application of friction damping to attenuate the vibrations of a complex structure was simplified to use the HBM by retaining one harmonic only. Later, the authors proved the applicability of the improved strategy to large Finite Elements (FE) models purposely reduced according to the industrial needs of fast and accurate calculation of the nonlinear forced response ([12],[16],[19]).

In this paper, the authors extend the improved strategy to the HBM with Multiple harmonics (MHBM) and apply it to the simplest lumped parameter system that is able to stress the dependence of the static and variable contact loads. The accuracy of the proposed coupled strategy is proved by using as a reference the solution obtained by means of the Direct Time Integration of the equation of motion.

Despite its simplicity, the analyzed system can be considered as the fundamental brick for the construction of highly detailed FE models of mechanical components with joint interfaces. As a consequence, the authors believe that the results obtained in this paper may help the analyst in the choice of the proper solution strategy to solve high detailed FE models or, at least, to be fully aware of the possible inaccuracy of their analysis when the uncoupled approach is chosen.

## 2  Governing equations of a system with joint interfaces

The governing equations of a vibrating system with joint interfaces in the time domain are:

$$\mathbf{m}\ddot{\mathbf{q}}(t) + \mathbf{c}\dot{\mathbf{q}}(t) + \mathbf{k}\mathbf{q}(t) = \mathbf{f}(t) + \mathbf{f}_C(\mathbf{q},\dot{\mathbf{q}},t) \qquad (1)$$

where $\mathbf{m}$, $\mathbf{c}$ and $\mathbf{k}$ are the mass, viscous linear damping and stiffness matrices, $\mathbf{q}$ the vector of the system degrees of freedom (dofs), $\mathbf{f}$ the external force and $\mathbf{f}_c$ the vector of the contact forces acting at the joint interfaces.

In case of a periodical excitation, the steady-state solution can be reached either by time integration of the balance equation (1) or by the Multi-Harmonic Balance Method (MHBM).
In the latter case, periodical quantities are expressed as a truncated series of harmonic terms as

$$\mathbf{q}(t) = \hat{\mathbf{q}}^{(0)} + \mathrm{Re}\left(\sum_{h=1}^{H} \hat{\mathbf{q}}^{(h)} e^{ih\omega t}\right)$$

$$\mathbf{f}(t) = \hat{\mathbf{f}}^{(0)} + \mathrm{Re}\left(\sum_{h=1}^{H} \hat{\mathbf{f}}^{(h)} e^{ih\omega t}\right) \qquad (2)$$

$$\mathbf{f}_c(\mathbf{q}, \dot{\mathbf{q}}, t) = \hat{\mathbf{f}}_c^{(0)}(\hat{\mathbf{q}}) + \mathrm{Re}\left(\sum_{h=1}^{H} \hat{\mathbf{f}}_c^{(h)}(\hat{\mathbf{q}}) e^{ih\omega t}\right)$$

and the set of non-linear differential balance equations (1) is turned in a set of non-linear algebraic equations with complex coefficients

$$\begin{aligned} \mathbf{k}\hat{\mathbf{q}}^{(0)} &= \hat{\mathbf{f}}^{(0)} - \hat{\mathbf{f}}_c^{(0)}(\hat{\mathbf{q}}) \\ \left(\mathbf{k} - h^2\omega^2 \mathbf{m} + ih\omega \mathbf{c}\right)\hat{\mathbf{q}}^{(h)} &= \hat{\mathbf{f}}^{(h)} - \hat{\mathbf{f}}_c^{(h)}(\hat{\mathbf{q}}) \end{aligned} \qquad \text{with} \quad h = 1 .. H \qquad (3)$$

where the two sets of equations represent the static and the dynamic governing equations of the system respectively, coupled to each other by the Fourier coefficients of the non-linear contact force $\mathbf{f}_c$, which depends on the Fourier coefficients of the displacement $\mathbf{q}$.

## 3   Contact forces

The solution of equation (2) requires the calculation of the contact forces $\mathbf{f}_c$ acting on the joint interface. To model the contact force, the state-of-the-art node-to-node contact element (Fig. 1) is used. The local contact stiffness is modelled by two linear springs $k_t$ and $k_n$ in the tangential and normal direction, respectively. In order to take into account friction phenomena occurring in the tangential direction, a slider is used to connect the bodies in contact in the tangential direction. According to the Coulomb friction law, when the tangential force $T(t)$ exceeds the limit value $\mu N(t)$, being $\mu$ the coefficient of friction, the slider moves with respect to the wall and the amount of slip between the contact nodes is $w(t)$. For each contact element, the tangential force T(t) and the normal force N(t) depend on the relative tangential and normal displacements of the contact points of the bodies in contact, named *u(t)* and *v(t)* respectively.

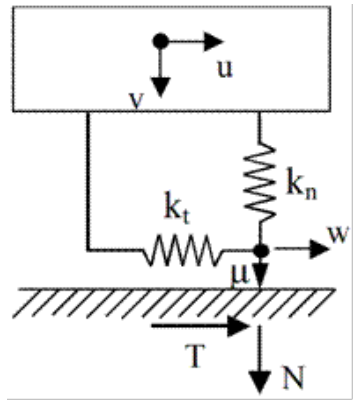

Fig. 1. Node-to-node contact element.

If the set of differential balance equations (1) are integrated in time, then at each time $t$, the normal contact force $N(t)$ is defined as

$$N(t) = \max(k_n v(t), 0) \tag{4}$$

where $k_n$ is the local normal stiffness of the contact. In case of negative values of $v(t)$, no traction forces are allowed but lift-off occurs.

A different approach is necessary along the tangential direction, since the value of the tangential contact force $T(t)$ depends on the contact state:

$$T(t) = \begin{cases} k_t[u(t) - w(t)] & \text{stick state} \\ \mu N(t) \text{sign}(\dot{w}) & \text{slip state} \\ 0 & \text{lift-off state} \end{cases} \tag{5}$$

where $k_t$ is the local tangential contact stiffness, $\mu$ is the coefficient of friction, $w(t)$ is the amount of actual slip. At each time $t$, a predictor step is performed ([8]), assuming the contact is in stick conditions

$$T^P(t) = k_t[u(t) - w(t)] = k_t[u(t) - w(t - \Delta t)] \tag{6}$$

where $\Delta t$ is the time step; then a corrector step is performed and the value of $T^P(t)$ is converted into the actual value $T(t)$ as:

$$T(t) = \begin{cases} T^P(t) & \text{stick state} \\ \mu N(t) \text{sign}(T^P(t)) & \text{slip state} \\ 0 & \text{lift-off state} \end{cases} \tag{7}$$

and the slider displacement $w(t)$ is accordingly computed as

$$w(t) = \begin{cases} w(t - \Delta t) & \text{stick state} \\ u(t) - \mu N(t) \text{sign}(T(t))/k_t & \text{slip state} \\ u(t) & \text{lift-off state} \end{cases} \tag{8}$$

In detail, the slip stare occurs if the modulus of the tangential force exceeds the Coulomb limit $\mu N(t)$, the lift-off state occurs if the normal force $N(t)$ is zero and the stick state occurs if none of the previous condition occurs.

Due to the use of a predictor-corrector strategy, an approximation is introduced in the calculation of the tangential contact force, since transitions between contact states and therefore the duration of each contact state are not exactly computed as in [2]. Nevertheless with a judicious choice of the time step Δt ( e.g. 100 time steps per cycle) the approximation is negligible for engineering purposes.

If the HBM is used and the set of algebraic equations (3) is used, the inputs are the Fourier coefficients of the relative displacements $\hat{u}$ and $\hat{v}$, computed by means of equation (13), while the outputs are the Fourier coefficients of the contact forces $\hat{\mathbf{f}}_c$.

The calculation procedure is an alternate time/frequency (AFT) method ([21],[22]), based on the following steps, as shown in Fig. 2:

1. Periodical relative displacements $u(t_j)$ and $v(t_j)$ are computed from Fourier coefficients $\hat{u}$ and $\hat{v}$ by Inverse Fast Fourier Transform (IFFT), being $t_j = j\Delta t$ with $j = 1...J$.
2. At each time step, contact forces $N(t_j)$ and $T(t_j)$ are computed by means of equations (4) - (8).
3. Fourier coefficients $\hat{N}$ and $\hat{T}$ are computed from $N(t_j)$ and $T(t_j)$ respectively by Fast Fourier Transform (FFT) and the vector of Fourier coefficients of contact forces $\hat{\mathbf{f}}_c$ is assembled.

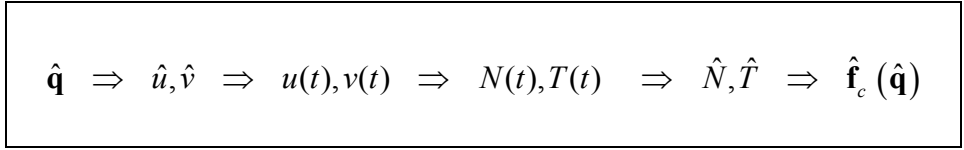

Fig. 2. MHBM, AFT method for the Fourier coefficients of contact forces.

## 4    Solution methods

In this section, the different solution methods, that can be used to find the periodical response of a vibrating system with joint interfaces under periodical excitation, are outlined.

The first method is the Direct Time Integration (DTI) of the balance equations. It gives the 'exact' value of the system response and it is usually used in the literature as a benchmark for approximate methods, which are developed in order to reduce the calculation time.

As already stated in Section 2, an approximate method largely used in the solution of dynamic problems with friction contact is the Multi-Harmonic Balance Method (MHBM). In the classical approach, the static balance equation is uncoupled from the dynamic balance equations; the set of equations (3) becomes

$$\mathbf{k}\,\hat{\mathbf{q}}^{(0)} = \hat{\mathbf{f}}^{(0)} - \hat{\mathbf{f}}_c^{(0)}\left(\hat{\mathbf{q}}^{(0)}\right)$$
$$\left(\mathbf{k} - h^2\omega^2\mathbf{m} + i h \omega \mathbf{c}\right)\hat{\mathbf{q}}^{(h)} = \hat{\mathbf{f}}^{(h)} - \hat{\mathbf{f}}_c^{(h)}\left(\hat{\mathbf{q}}^{(0)},\hat{\mathbf{q}}^{(h)}\right) \quad \text{with } h = 1 .. H \quad (9)$$

The static balance equation is solved and then the 0th order coefficients of the displacements $\hat{\mathbf{q}}^{(0)}$ are used as inputs for the solution of the dynamic equations. As a consequence, the value of the static contact forces, computed in the static balance equations is no more updated when the dynamic balance is imposed.

In a more general approach, the static and the dynamic balance equations formulated by the MHBM are solved together and the overall balance conditions (static + dynamic) are computed by means of equations

$$\begin{aligned} \mathbf{k}\,\hat{\mathbf{q}}^{(0)} &= \hat{\mathbf{f}}^{(0)} - \hat{\mathbf{f}}_c^{(0)}\!\left(\hat{\mathbf{q}}^{(0)}, \hat{\mathbf{q}}^{(h)}\right) \\ \left(\mathbf{k} - h^2\omega^2 \mathbf{m} + i\,h\,\omega\,\mathbf{c}\right)\hat{\mathbf{q}}^{(h)} &= \hat{\mathbf{f}}^{(h)} - \hat{\mathbf{f}}_c^{(h)}\!\left(\hat{\mathbf{q}}^{(0)}, \hat{\mathbf{q}}^{(h)}\right) \end{aligned} \quad \text{with} \quad h = 1\,..\,H \qquad (10)$$

which correspond to the fully coupled equation set (3). This procedure will be referred as the 'coupled' calculation of the forced response, where the term coupled means that the static and the dynamic quantities involved in the dynamic response of the vibrating system are updated at each iteration and may change for different excitation frequencies.

In the following section, it is shown by means of a set of numerical examples, that the calculation of the forced response by using an uncoupled frequency domain based procedure which assumes the values of the static quantities (forces and displacements) not affected by the nonlinear contact model represented by equations (4)-(8) may lead to uncertainties and inaccuracies in the prediction of the maximum amplitude response.

## 5    Application

The uncoupled and the coupled HBM approaches are here applied to a 2 dofs model, shown in Fig. 3, made of a mass $m$ with 2 degrees of freedom (horizontal and vertical displacements $q_x(t)$ and $q_y(t)$). The mass is connected to ground by means of two springs ($k_x$ and $k_y$) and two viscous damping elements ($c_x$ and $c_y$). The reader can refer to Table 1 for the complete list of the system parameters. Despite its simplicity, the proposed test case can be considered as the fundamental brick for the construction of highly detailed FE models of mechanical components with joint interfaces ([3]-[8],[10]-[12]).

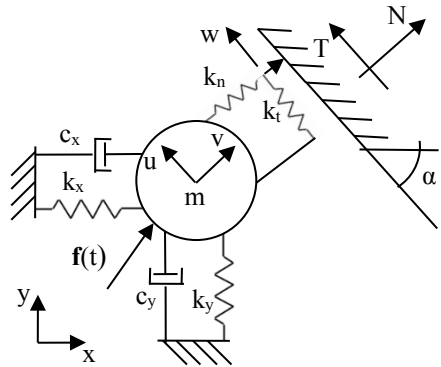

Fig. 3. Test case: 2 dofs model.

Table 1 : lumped parameters system

| Parameter | value |
|---|---|
| $m$ | 1 kg |
| $k_x$ | 3.9e5 N m$^{-1}$ |
| $k_y$ | 2.5e6 N m$^{-1}$ |
| $c_x$ | 62.8 N (m/s)$^{-1}$ |
| $c_y$ | 62.8 N (m/s)$^{-1}$ |
| $k_t$ | 3e5 N m$^{-1}$ |
| $k_n$ | 3e5 N m$^{-1}$ |
| μ (coefficient of fiction) | 0.5 |
| $F_{pl}$ | 420 N |
| $|F_{exc}|$ | 1, 4, 6, 18, 24, 30, 42, 60 N |

An external force vector **f**(*t*) is applied to the mass and during the motion the mass may enter in contact with the fixed wall whose slope is defined by the angle *α* located on its right-hand side. A node-to-node contact element is used to model the periodical contact forces.

The matrices and vectors in the governing equations are

$$\mathbf{m}\ddot{\mathbf{q}}(t) + \mathbf{c}\dot{\mathbf{q}}(t) + \mathbf{k}\mathbf{q}(t) = \mathbf{f}(t) - \mathbf{f}_c(\mathbf{q},\dot{\mathbf{q}},t)$$

$$\mathbf{m} = \begin{bmatrix} m & 0 \\ 0 & m \end{bmatrix}; \mathbf{c} = \begin{bmatrix} c_x & 0 \\ 0 & c_y \end{bmatrix}; \mathbf{k} = \begin{bmatrix} k_x & 0 \\ 0 & k_y \end{bmatrix}; \mathbf{q} = \begin{Bmatrix} q_x \\ q_y \end{Bmatrix}; \mathbf{f} = \begin{Bmatrix} f_x \\ f_y \end{Bmatrix}; \mathbf{f}_c = \begin{Bmatrix} f_{cx} \\ f_{cy} \end{Bmatrix} \quad (11)$$

where the vector of contact forces $\mathbf{f}_c(\mathbf{q},\dot{\mathbf{q}},t)$ is defined as

$$\mathbf{f}_c(\mathbf{q},\dot{\mathbf{q}},t) = \begin{Bmatrix} f_{cx} \\ f_{cy} \end{Bmatrix} = \begin{Bmatrix} N(\mathbf{q},\dot{\mathbf{q}},t)\sin(a) - T(\mathbf{q},\dot{\mathbf{q}},t)\cos(a) \\ N(\mathbf{q},\dot{\mathbf{q}},t)\cos(a) + T(\mathbf{q},\dot{\mathbf{q}},t)\sin(a) \end{Bmatrix} \quad (12)$$

and the relative displacements of the mass with respect to the fixed wall are defined as:

$$\begin{Bmatrix} u \\ v \end{Bmatrix} = \begin{Bmatrix} -q_x \cos(a) + q_y \sin(a) \\ q_x \cos(a) + q_y \sin(a) \end{Bmatrix} \quad (13)$$

The 2-dof model of Fig. 3 is used in two different configurations to gradually increase the complexity of the system response. In the first case a consistency check of the two methods in the frequency domain is performed by comparing the results obtained by the coupled and the uncoupled approach with DTI. Moreover, the typical behaviour of a vibrating structure in presence of a sliding contact is presented. In the second case the main limitations and drawbacks caused by the uncoupled approach in the frequency domain are highlighted and it is proved that the coupled approach gives the exact solution for a vibrating system with one friction contact. In the first configuration, the fixed wall is horizontal ($\alpha=0°$, Fig. 4). The natural frequencies of the linear system are 100 Hz and 250 Hz, corresponding to two orthogonal mode shapes $\varphi_1=[A_x\ 0]$ and $\varphi_2=[0\ A_y]$. A constant pre-load $F_{pl}$ is applied to the mass $m$ along the direction normal to the contact (y-axis) and an exciting force $F_{exc}(t)$ varying harmonically with a given frequency $f_{exc}$ is applied to $m$ along the direction tangential to the contact (x-axis). All the loads applied to the system (external and contact loads) are aligned with the two mode shapes $\varphi_1$ and $\varphi_2$, therefore there is no coupling between the x and y motion. As a consequence, this simple case allows for a constant normal load during vibrations.

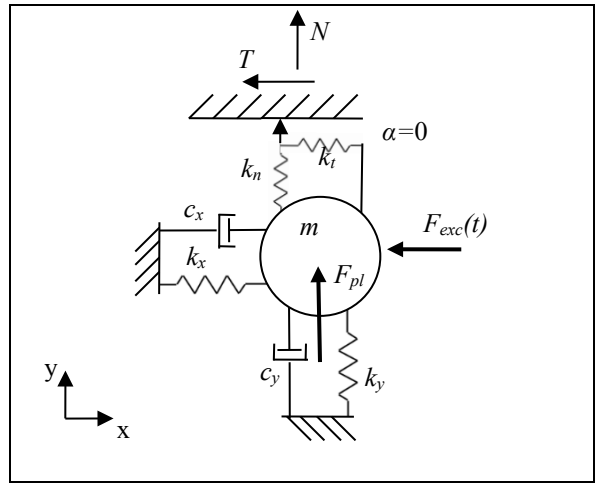

Fig. 4. Case of a horizontal fixed wall ($\alpha=0°$).

The pseudo Frequency Response Functions (pseudo-FRFs) of the system calculated as the ratio between the amplitude of the multi-harmonic response and the excitation amplitude by means of the classical uncoupled (U) and coupled (C) MHBM approaches are plotted together at 8 different excitation force amplitudes in Fig. 5. When the uncoupled approach is used the governing equation (9).a is previously solved by loading the system with the static preload $\hat{\mathbf{f}}^{(0)} = \begin{bmatrix} 0 & F_{pl} \end{bmatrix}^T$ in order to calculate the static contact forces $\hat{\mathbf{f}}_c^{(0)} = \begin{bmatrix} T^{(0)} & N^{(0)} \end{bmatrix}^T = \begin{bmatrix} 0 & 46\ \text{N} \end{bmatrix}^T$ and the static displacements $\hat{\mathbf{q}}^{(0)} = \begin{bmatrix} 0 & y^{(0)} \end{bmatrix}^T = \begin{bmatrix} 0 & 0.15\text{e}-3\ \text{m} \end{bmatrix}^T$ which are used in equation (9).b to calculate the dynamic response. A number of harmonics H=5 is retained for both the frequency domain based approaches. It is possible to see that in this case the two approaches give the same results (U=C) for each comparison: by increasing the amplitude of the excitation force both the value of the peak amplitude and the corresponding frequency change. In detail, when the force is small with respect to the normal preload ($|F_{exc}|$=1 N; 4 N vs. $F_{pl}$=420 N), the contact is in fully stick condition and equation (5).a holds during the whole vibration as it is possible to see in Fig. 6.a where the tangential force $T(t)$ is plotted versus the relative displacement $u(t)$ for the peak frequency ($f_{exc}$=132 Hz). In this case $T(t)$ is proportional to $u(t)$, therefore the calculation is linear and a global stiffening effect is observed with respect to the natural frequency of the system without contact (free). The reason is due to the presence of the stiffness $k_t$ which is added in parallel to the structural stiffness $k_x$. By increasing $|F_{exc}|$ the peak amplitude decreases to a minimum ($|F_{exc}|$=24 N) since the amount of relative tangential motion $u(t)$ at the contact is large enough to let the contact slip (see Fig. 6.b where the hysteresis loop at the contact is plotted for $|F_{exc}|$=24 N and $f_{exc}$=115 Hz). The peak amplitude increases for higher values of $|F_{exc}|$ ($|F_{exc}|$=30, 42, 60 N) since the damping generated by friction is no more sufficient to limit the structure vibrations. Fig. 6.c showing the hysteresis loop for $|F_{exc}|$=60 N and $f_{exc}$=101 Hz can better explain the result if compared with the hysteresis loop of Fig. 6.b: the increase of the tangential displacements $u(t)$ is higher than the increase of the dissipated energy, i.e. the area of the hysteresis loop.

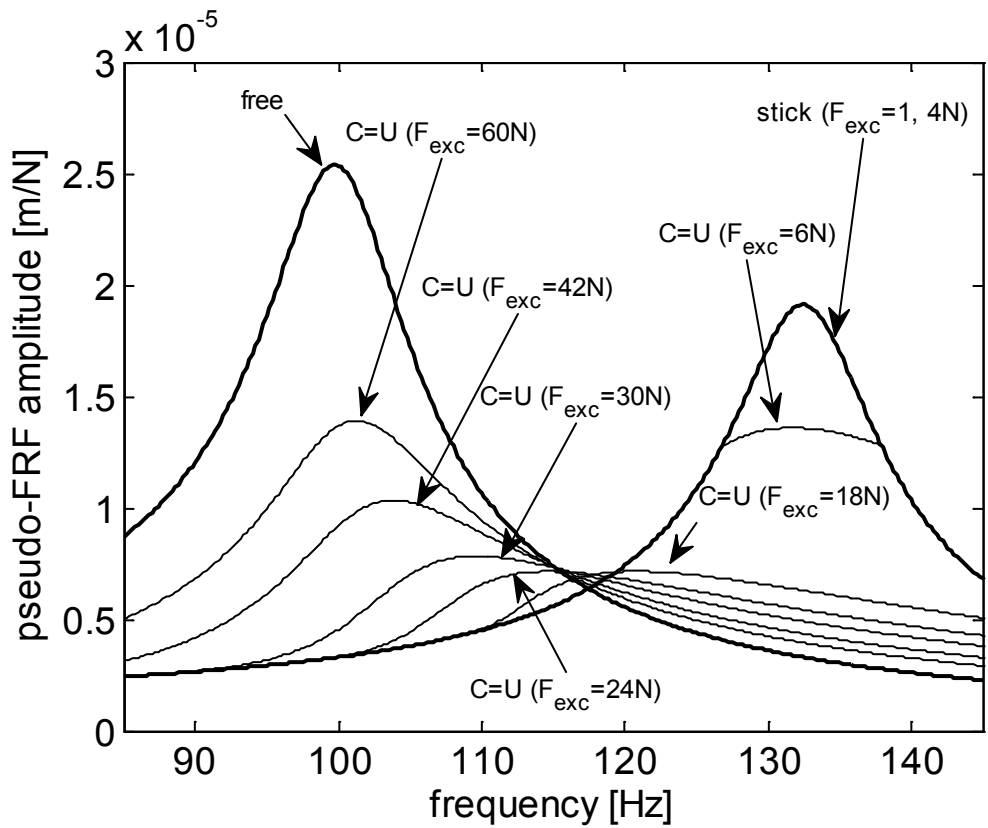

Fig. 5. $\alpha$= 0°, pseudo-FRFs for different amplitude of $F_{exc}$.

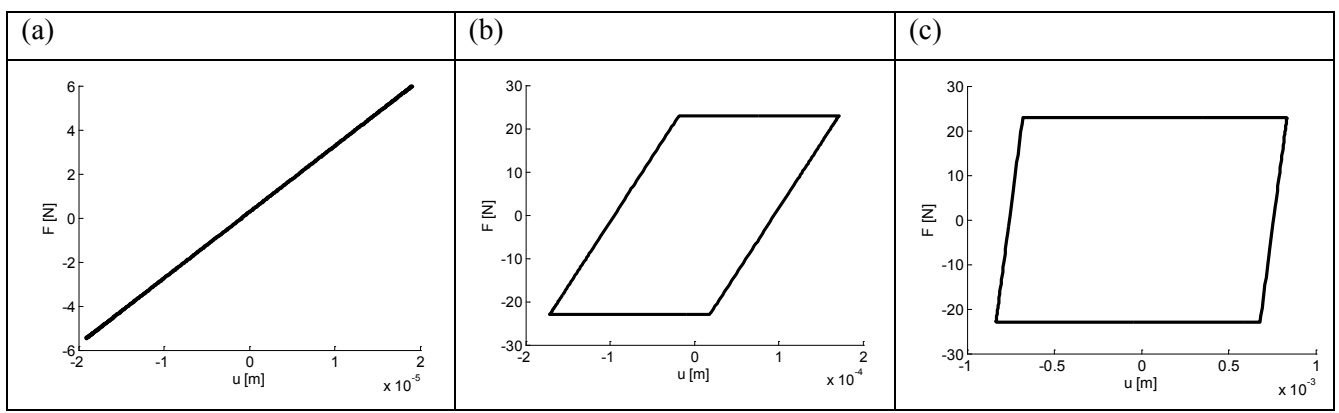

Fig. 6. $\alpha$= 0°, hysteresis loops for three different resonant responses: (a) $|F_{exc}|$=1 N; (b) $|F_{exc}|$=24 N; (c) $|F_{exc}|$=60 N.

The direct time integration (DTI) is used as a consistency check of the two methods U and C. In detail, the calculation is performed according to two load steps: in the first load step only the preload is gradually applied through a quasi-static ramp until a constant value $F_{pl}$ is reached. In the second load step the harmonic exciting force is gradually applied in order to reach the amplitudes as listed in Table I and at the corresponding peak frequency. The perfect matching of the three methods (U, C, DTI) can be globally shown in terms of optimization curve (Fig. 7) where the peak amplitude of each FRF is plotted versus the $F_{pl}/F_{exc}$ ratio. The solid line connects the maxima obtained with the two frequency domain approaches while circle markers correspond to the amplitude of oscillation of the DTI stationary solutions divided by $|F_{exc}|$. Note the transition from full stick to stick-slip state for $F_{pl}/F_{exc}= 106$.

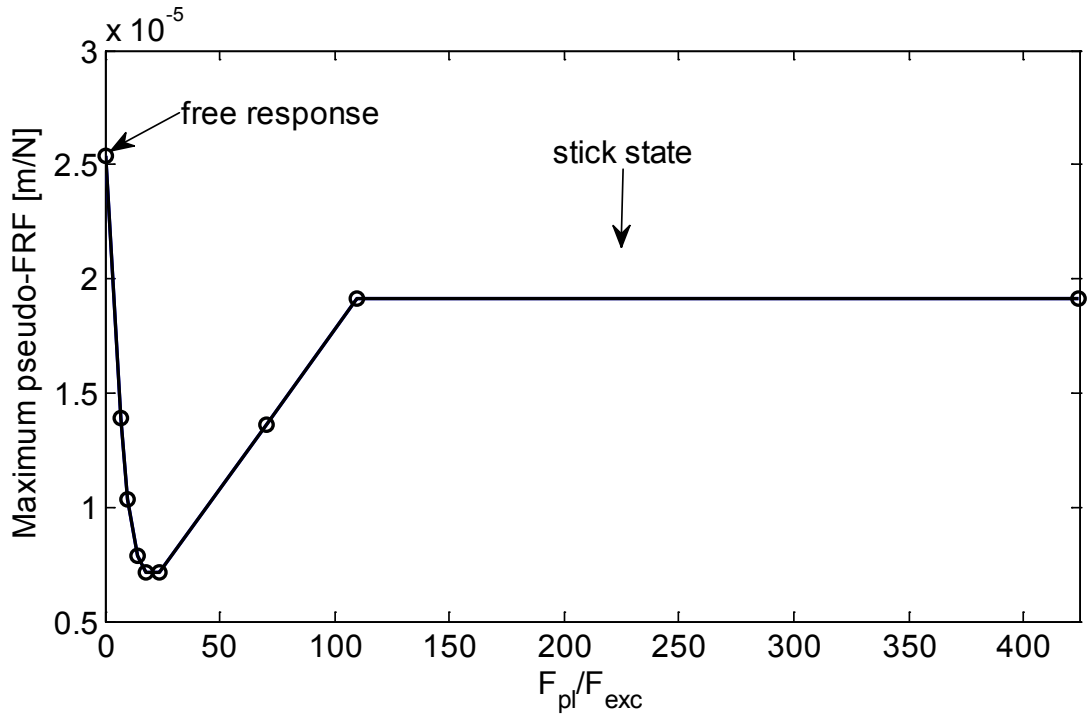

Fig. 7. Optimization curve ($\alpha= 0°$)

The second configuration that is presented is this paper is a more general case where the static preload $F_{pl}$ is not normal to the contact and both the variable external force $F_{exc}(t)$ and the contact forces $T(t)$ and $N(t)$ are applied along directions which couple the x and y motion. For this reason $\alpha= 45°$ and the preload $F_{pl}$ acts again along the y-axis (Fig. 8. General case of contact interaction ($\alpha= 45°$).). In this case the relationship between $F_{pl}$ and the static contact forces (and static displacements) is not unique as in the previous case and it may depend on the history of application of $F_{pl}$. This is an important issue to discuss in order to use the uncoupled approach which needs $\hat{\mathbf{q}}^{(0)}$ to solve the set of equations (9).b.

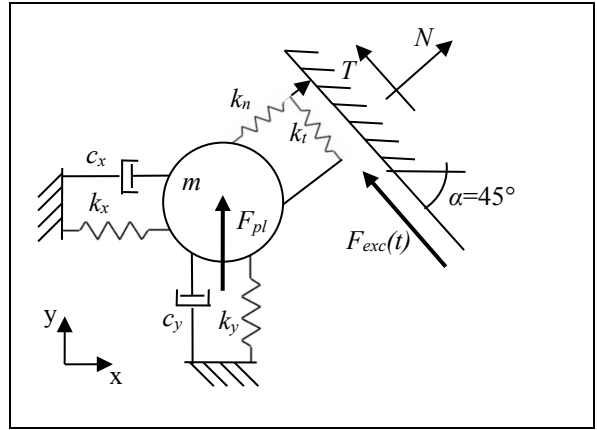

Fig. 8. General case of contact interaction ($\alpha= 45°$).

As an example, two sub-cases (U1 and U2) are analyzed where the DTI integration is used to calculate the static quantities $\hat{\mathbf{q}}^{(0)}$ and $\hat{\mathbf{f}}_c^{(0)}$. In sub-case U1 (Fig. 9.a) the static preload $F_{pl}$ is applied with a monotonic increasing ramp until the final value is reached ($F_{pl}=$ 420 N) and the transient response has terminated, therefore the static contact forces $\hat{\mathbf{f}}_{c,1}^{(0)} = \begin{bmatrix} T_1^{(0)} & N_1^{(0)} \end{bmatrix}^T$ (Fig. 9.b) and the static displacements $\hat{\mathbf{q}}_1^{(0)} = \begin{bmatrix} x_1^{(0)} & y_1^{(0)} \end{bmatrix}^T$ (Fig. 9.c) are stored in order to solve the set of equations (9).b. It must be noted that subcase U1 corresponds to solve the first set of equation 11.a since the preload is applied with a single slope ramp. In this case the preload is large enough to determine a positive slip of the contact ($T^{(0)}=\mu N^{(0)}$). In the sub-case U2 (Fig. 9.d) the static preload is applied with a different history profile in order to obtain at the end of the load step the same value ($F_{pl}=$ 420 N) but two new values of the static contact forces $\hat{\mathbf{f}}_{c,2}^{(0)} = \begin{bmatrix} T_2^{(0)} & N_2^{(0)} \end{bmatrix}^T$ (Fig. 9.e) and static displacements $\hat{\mathbf{q}}_2^{(0)} = \begin{bmatrix} x_2^{(0)} & y_2^{(0)} \end{bmatrix}^T$ (Fig. 9.f). In particular the decreasing ramp determines a negative slip of the contact. Now the two sets $\hat{\mathbf{q}}_1^{(0)}$ and $\hat{\mathbf{q}}_2^{(0)}$ are used to solve equations (9).b according to the uncoupled approach and the FRFs are compared with the FRFs obtained by the coupled approach which does not require the static displacements as input parameters. The comparison is shown in terms of optimization curves where it is possible to see that, for given values of $F_{pl}/F_{exc}$ ratios, the scatter between sub-case U1 and U2 is large and the predicted peak amplitudes are different from the unique solution calculated by the coupled approach C. In particular, the optimization curve corresponding to subcase U1 shows higher peak responses than the optimization curve of subcase U2. This can be explained since the static normal load of U1 is higher than the static normal load of U2, therefore for subcase U1 the contact enters the slip state for excitation forces higher than those of subcase U2.

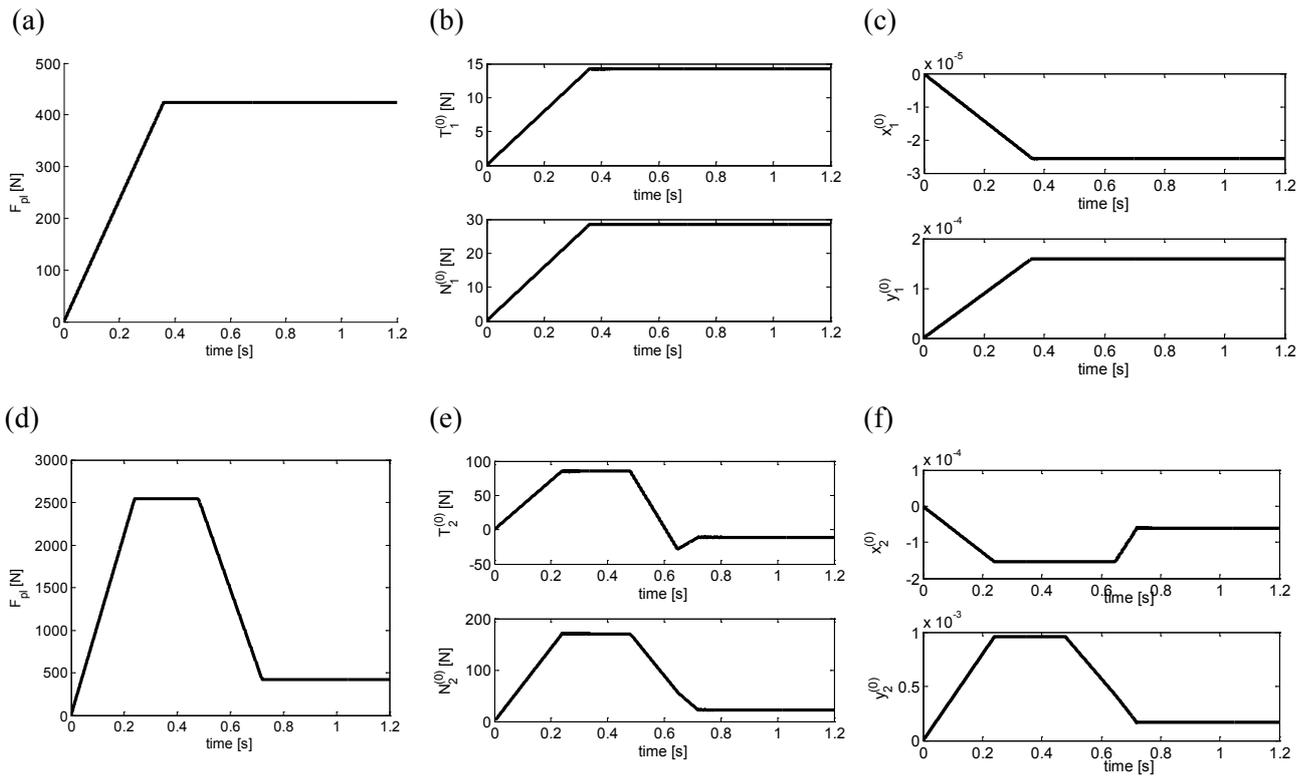

Fig. 9. Direct time integration, application of the preload $F_{pl}$ according to two different load profiles: (a) first preload profile 'U1'; (b) contact loads evolution U1 (final values $T_1^{(0)}$=14.2 N; $N_1^{(0)}$=28.5 N) ; (c) displacements evolution U1 (final values $x_1^{(0)}$=-2.6e-5 m ; $y_1^{(0)}$=1.6e-4 m); (d) second preload profile 'U2'; (e) contact loads evolution U2 (final values $T_2^{(0)}$=-10.8 N; $N_2^{(0)}$=22.9 N); (f) displacements evolution U2 (final values $x_2^{(0)}$=-6e-5 m ; $y_2^{(0)}$=1.7e-4 m).

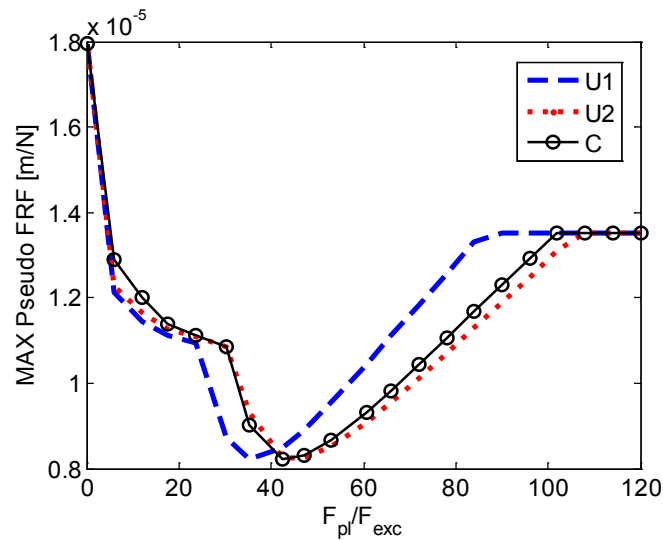

Fig. 10. Optimization curves with different approaches: Uncoupled approach (U1 and U2), coupled approach (C), direct time integration (circle markers).

The FRFs and the corresponding hysteresis loops obtained at the frequency where the peak responses occur are now plotted for $F_{pl}/F_{exc}$= 84 and 30 in Fig. 11.a and Fig. 11.b (FRFs) and Fig. 12.a and Fig. 12.b (hysteresis loops). In the first case ($F_{pl}/F_{exc}$= 84) the peak amplitude for the subcase U1 corresponds to a fully stick contact as it is clearly shown by the corresponding linear dependence of the tangential force $T(t)$ from the tangential displacement $u(t)$ in Fig. 12.a. The peak response of subcase U2 is similar to the peak response calculated by the coupled approach C since the dissipated energy is similar for the two cases as it is visible in Fig. 12.a. Nonetheless, the static tangential displacements is correctly calculated by the coupled approach only since the hysteresis loop perfectly match with the hysteresis loop calculated by the DTI. Moreover, it must be noted that the hysteresis loops do not keep the typical shape of Fig. 6.b and Fig. 6.c since the normal load $N(t)$ varies with time according to equation (4). In the second case ($F_{pl}/F_{exc}$= 30, Fig. 11.b) the peak amplitude of subcase U1 markedly differs from the peak amplitude obtained with C and subcase U2. The reason is found in Fig. 12.b where the hysteresis loop calculated at the peak frequency for U1 shows a full contact during vibration while, on the contrary, subcase U2 and C shows partial lift-off during vibration which decreases the damping effectiveness of the contact producing higher response amplitudes. It must be noted, however, that only the coupled approach predict exactly both the static and the dynamic response of the system calculated by the DTI as well as the corresponding hysteresis loop.

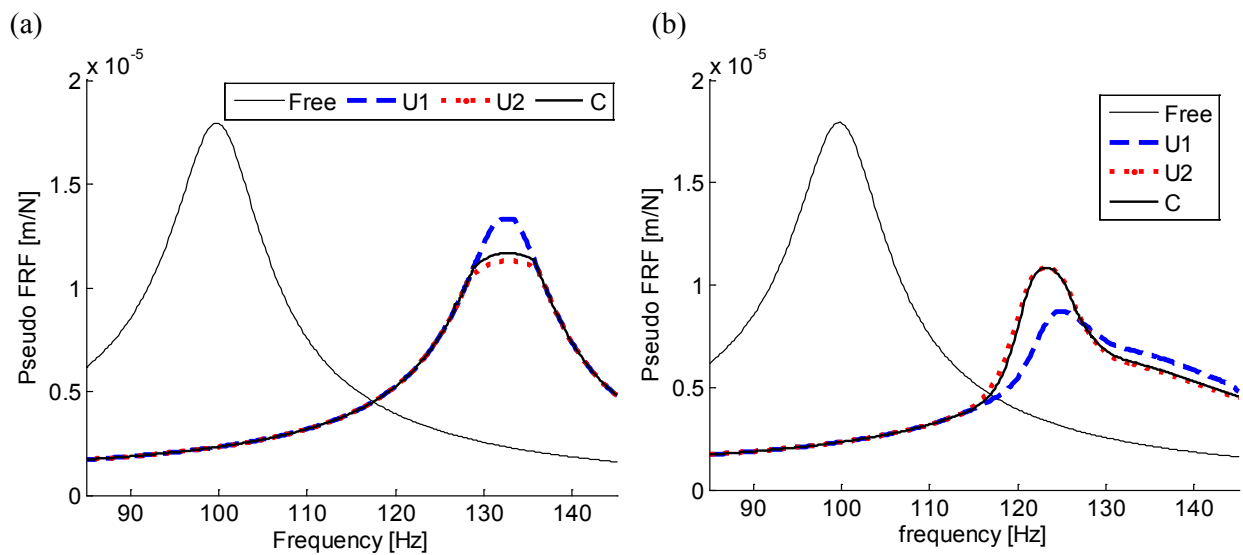

Fig. 11. Pseudo-FRFs comparison for different $F_{pl}/F_{exc}$ ratios: (a) $F_{pl}/F_{exc}$=84, (b) $F_{pl}/F_{exc}$=30.

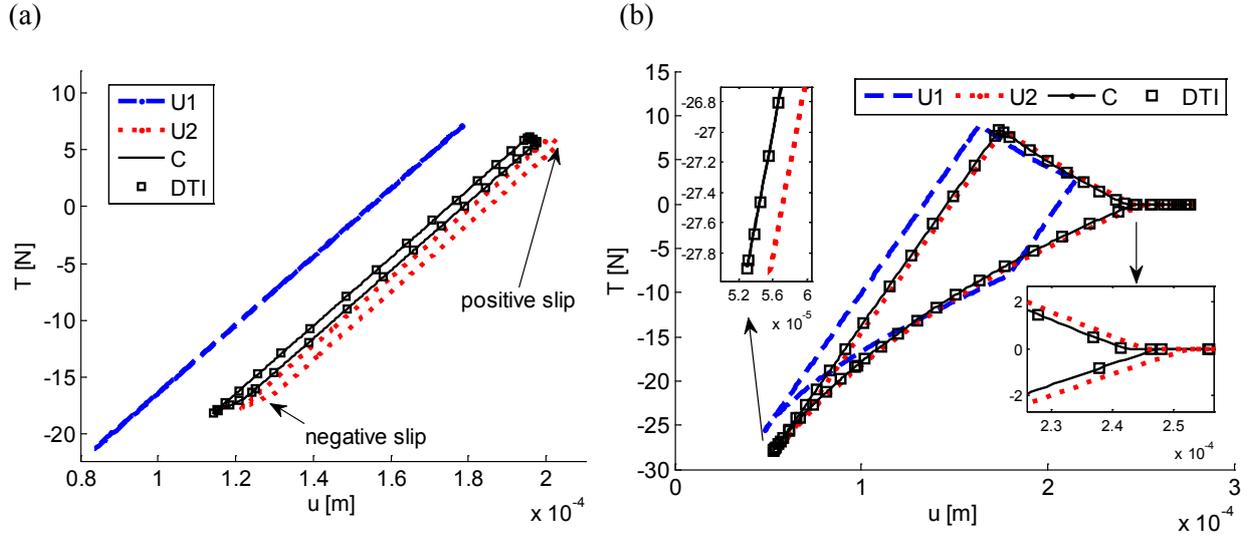

Fig. 12. Hysteresis cycles for the peak frequency U1, U2, C and comparison with DTI: (a) $F_{pl}/F_{exc}$=84, (b) $F_{pl}/F_{exc}$=30.

As a final result, the uniqueness of the solution generated by the coupled approach is verified in one case ($F_{pl}/F_{exc}$=84) with the DTI, by applying the harmonic load $F_{exc}(t)$ as a second load step following the load steps of Fig. 9b and Fig. 9c, and by calculating the steady state response of the system for the two cases. The static preload $F_{pl}$ is kept constant to 420 N during all the simulation. The excitation frequency is equal to the peak frequency ($f_{exc}$= 132 Hz). The transient response of the system is plotted in terms of contact loads in Fig. 13.a and in terms of displacements in Fig. 13.b. It is possible to see that the system starts to vibrate with different initial conditions, in particular the red curves refer to the initial condition of subcase U1 while black curves refer to the initial condition of subcase U2. After the transient response has terminated, the two solutions converge to one solution that is also the unique solution found by the coupled approach in the frequency domain. It is possible to note that the final static values of the contact forces $T^{(0)}$ and $N^{(0)}$ are different from both the static contact forces at the end of the first load step of figure 8.b ($T_1^{(0)}, N_1^{(0)}$) and 8.e ($T_2^{(0)}, N_2^{(0)}$). The final correct values depend on the strong nonlinear coupling of the static and dynamic quantities determined by the friction contact and cannot be calculated separately from the dynamic calculation, otherwise large discrepancies may be found in the final results with respect to the correct solution.

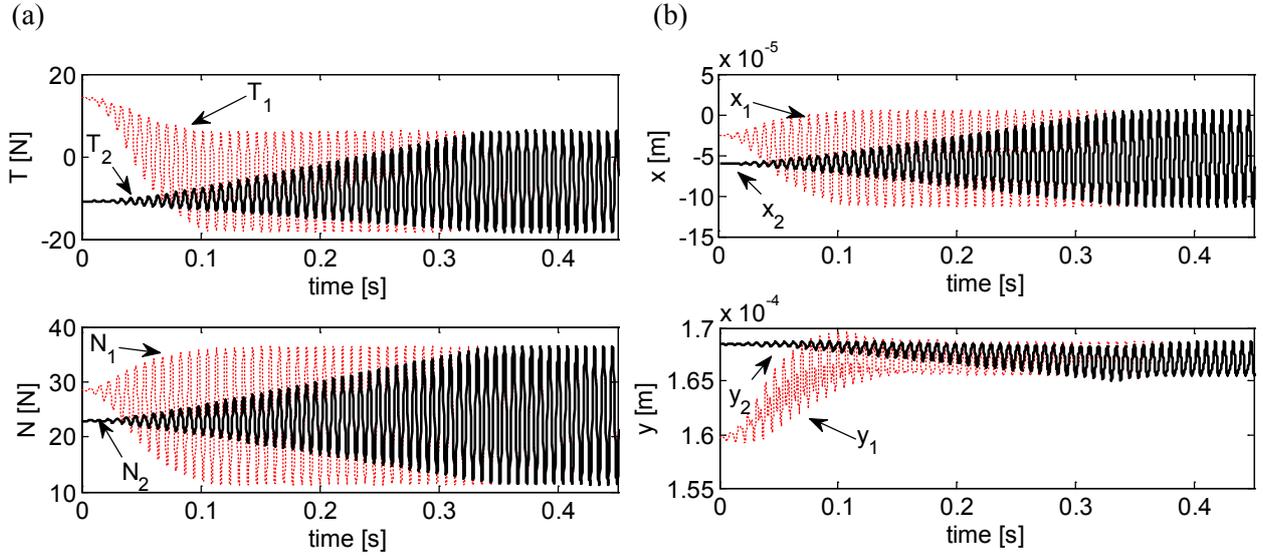

Fig. 13. Transient response for different initial conditions (U1, dotted curves and U2, solid curves) and common steady state response of the system after the application of a harmonic excitation force: (a) contact forces (steady state: $T^{(0)}$= -6 N; $N^{(0)}$= 23.9 N), (b) displacements (steady state: $x^{(0)}$= -5.4e-5 m; $y^{(0)}$= 1.7e-4 m).

## 6 Conclusions

In this paper a numerical investigation is presented to analyze and compare uncoupled and coupled static and dynamic approaches to calculate the nonlinear forced response of structures with friction contacts by means of the MHBM, while the DTI is used as a reference for the MHBM approaches.

According to the uncoupled approach, the static governing equations are solved in advance in order to compute the static tangential and normal contact forces, that acts on the contact interfaces, and then dynamic governing equations are solved by by keeping the static contact forces fixed. is based on the preliminary solution of the static governing equations of the system and then on the solution of the dynamic governing equations, by using the outputs of the static analysis as input parameters.

The coupled approach, on the contrary, is based on the simultaneous solution of both the static and the dynamic governing equations.

A simple lumped parameter system has been used to compare the two approaches. It has been shown that:

- If the uncoupled approach is used, different load step sequences in the preliminary static analysis may lead to different static solutions and consequently to non-unique dynamic solutions.
- If the coupled approach is used, whatever the initial condition is in terms of static contact loads and static displacements, the steady-state solution is unique when the slip state occurs during vibration.
- When no slip occurs, the system is linear and problem is trivial.

- If the tangential and the normal direction of the contact are uncoupled, the two approaches predict the same vibration amplitude, coincident with that of the DTI method.
- If the tangential direction of the contact is coupled to the normal direction, only the coupled approach leads to a solution perfectly coincident with that of the DTI method.
- The coupled approach avoids the application of two load steps during analysis (first static, then dynamics) since static displacements are an output of the method as well as the dynamic displacements.

The authors believe that the results obtained in this paper, although referred to a simple vibrating system, give useful insights in the problem of modelling joint interfaces in vibrating system and may help the analyst in the choice of the proper solution strategy to solve high detailed FE models or, at least, to be fully aware of the possible inaccuracy of their analysis when the uncoupled approach is chosen.